\documentclass[11pt,a4paper,leqno]{amsart}

\usepackage{amsmath,amssymb, amsfonts}

\def\<#1,#2>{\langle\,#1,\,#2\,\rangle}

\def\qed{\ensuremath{\hfill\Box}}
\def\qed{\ensuremath{\hfill\Box}}

\def\phi{\varphi}

\def\cala{\mathcal{A}}
\def\al{\alpha}
\def\eps{\varepsilon}

\theoremstyle{plain}
\newtheorem{Lemma}{Lemma}

\newtheorem{Theorem}[Lemma]{Theorem}

\theoremstyle{definition}

\theoremstyle{remark}
\newtheorem{Remark}[Lemma]{Remark}

\title{A property of the bidimensional sphere}
\author{Marius Cavachi}

\address{"Ovidius" University of Constan\c ta, 124 Mamaia Bvd
900527 Constan\c ta, Romania}
\email{mcavachi@yahoo.com}

\begin{document}

\begin{abstract}

It is natural to ask for a reasonable constant $k$ having the property that any open set of area greater than $k$ on a bidimensional sphere of area $1$ always contains the vertices of a regular tetrahedron. We shall prove that it is sufficient to take $k=\frac 34$. In fact we shall prove a more general result. The interested reader will not have any problem in establishing that $\frac 34$ is the best constant with this property.

\medskip

\end{abstract}

\maketitle

Our result is the following:

\begin{Theorem}
Let $n$ be a positive integer, and let $S$ be a bidimensional sphere of area 1. If $M\subset S$ is an open set
of area greater than $\displaystyle\frac{n-1}{n}$ and $X\subset S$ is a finite set with $n$ elements,  then there exists a rotation $\rho$ of the sphere such that
$\rho(X)\subset M$.
\end{Theorem}

In the proof, we use the following result whose proof we postpone:

\begin{Lemma}
Let $M, M'\subset S$ be open sets such that $\cala(M)>\cala(M')$\footnote{For any $A\subset S$,  $\cala(A)$ denotes its area.}. Then there exists a finite number of mutually disjoint spherical caps $U_\al$ and rotations $\rho_\al$ such that:

$(i)$ $\bigcup_\al U_\al\subset M$;

$(ii)$ $M'\subset \bigcup_\al \rho_\al(U_\al)$;

$(iii)$ $M\setminus \bigcup_\al U_\al$ has non-empty interior.
\end{Lemma}

\smallskip

\noindent\textbf{Proof of the Theorem.}~ Let $\mu$ be a Haar measure on $SO(3)$ such that $\mu(SO(3))=1$.

For any $A\subset S$, let $\Phi_A$ be the characteristic function of $A$.

Fix $a\in S$ and let $I_a^A\in \mathbb{R}$ be $I_a^A=\int_{SO(3)}\Phi_A\circ x(a)d\mu(x)$. 

\smallskip

\begin{Remark}
Note that if $b$ is an arbitrary point on $S$, then $I_a^A=I^A_b$. Indeed if $\rho\in SO(3)$ is such that $\rho(a)=b$ (and such a $\rho$ always exists), then:
\begin{equation*}
\begin{split}
I_b^A&=\int_{SO(3)}\Phi_A\circ x(\rho(a))d\mu(x)=\int_{SO(3)}\Phi_A\circ (x\circ \rho)(a))d\mu(x\circ \rho)\\
&=\int_{SO(3)}\Phi_A\circ x(a)d\mu(x),
\end{split}
\end{equation*}
since $d\mu(x\circ \rho)=d\mu(x)$, the Haar measure being rotation invariant.

Moreover, if $B\subset S$ is an open set such that there exists $\rho_1\in SO(3)$ with $\rho_1(A)=B$, then again $I_a^A=I^B_a$.
Indeed, 
\begin{equation*}
\begin{split}
I_a^B&=\int_{SO(3)}\Phi_{\rho(A)}\circ x(a)d\mu(x)=\int_{SO(3)}\Phi_A\circ \rho_1^{-1}\circ x(a))d\mu(x)\\
&=\int_{SO(3)}\Phi_A\circ (\rho_1^{-1}\circ x)(a)d\mu(\rho_1^{-1}\circ x)=I_a^A.
\end{split}
\end{equation*}
\end{Remark}

\smallskip

\noindent Returning to the problem, if $X=\{a_1,\ldots, a_n\}$, let 
$$f:SO(3)\rightarrow \mathbb{R},\qquad f(x)=\sum_{i=1}^n\Phi_M\circ x(a_i).$$ 
Note that it is enough to find an $x\in SO(3)$ with $f(x)> n-1$. Then, since $f(x)$ is an integer $\leq n$, we obtain  $f(x)=n$ and hence 
$x(a_1),\ldots, x(a_n)\in M$, which proves the Theorem. To find such an $x$, it is enough to show that $$\int_{SO(3)}f(x)d\mu(x)>n-1.$$

\noindent But this means that $$\sum_{i=1}^n I_{a_i}^M>n-1,$$ which is implied by $$I_{a_i}^M>\frac{n-1}{n}$$ for each $i$, that is $$I_{a}^M>\frac{n-1}{n}.$$

We divide the sphere $S$ in $n$ spherical lunes $F_1,\ldots, F_n$ of equal areas. Obviously, each $F_i$ can be  obtained as a rotation of $F_1$. This implies:
$$1=I_a^S=\sum _{i=1}^nI_a^{F_i}=nI_a^{F_1}, \quad \text{hence}\quad I_a^{F_1}=\frac 1n.$$
Let now $M'=S\setminus F_n$. Then 
$$I_a^{M'}=\sum _{i=1}^{n-1}I_a^{F_i}=\frac{n-1}{n}.$$
With $U_\al$ and $\rho_\al$ as in the Lemma, we deduce:
$$I_a^M>I_a^{\cup_\al U_\al}=\sum_\al I_a^{U_\al}=\sum_\al I_a^{\rho_\al(U_\al)}\geq I_a^{M'}=\frac{n-1}{n},$$
and the proof is complete.\qed

\vskip.5cm

\noindent{\bf Proof of the Lemma.}~ Let $0<m<1$ and let $C_i$, for $i\in\{1,\ldots, k\},$ be spherical caps of diameter $d$ such that $$\bigcup _{i=1}^kC_i=S,$$ and let $P_i$ be the plane containing the center of $S$ and parallel to the circle bounding $C_i$. If $\pi_i:S\rightarrow P_i$ is the orthogonal projection on $P_i$, we can choose $d$ small enough such that: 
\begin{itemize}
\item For any open $C\subset C_i$, we have $\cala(\pi_i(C))>m\cala(C).$
\item For any $A\neq B\in C_i$, we have the inequality of segment lengths: $$|\pi_i(A)\pi_i(B)|>m\cdot|AB|.$$
\end{itemize}

\noindent Define now 
$\displaystyle M_1=C_1\cap M$, $M_2=C_2\cap (M\setminus M_1)$, $M_3=C_3\cap (M\setminus M_1\cup M_2),\ldots$, $M_k=C_k\cap (M\setminus M_1\cup\cdots \cup M_{k-1}),$
and similarly construct $M'_1, M'_2,\ldots, M'_k$.

Let $N_i=\pi_i(M_i)$, $N'_i=\pi_i(M'_i)$. For $1-m$ close enough to $0$, we have:
$$\sum_{i=1}^k\cala(N_i)>\sum_{i=1}^k\cala(N'_i).$$
In each plane $P_i$, we fix a side length $\eps$ square lattice. It can be proven (see \cite[pag. 315,327]{M}) that the number $n_i$ of squares contained in $N_i$ is $$\displaystyle \frac{1}{\eps^2}\cala(N_i)+O(\frac 1\eps),$$ and analogously we have an approximation for the number $n'_i$ of  squares contained in $N'_i$. Hence, for small enough $\eps$, we get $$\sum_{i=1}^kn_i>\sum_{i=1}^kn'_i.$$ 

\noindent Therefore, we can choose an injection $u$ from the set $\mathcal{P}'$ of squares contained in $\bigcup_{i=1}^kN'_i$ into the set $\mathcal{P}$ of squares contained in $\bigcup_{i=1}^kN_i$. 

Let $P\in \mathcal{P}'$ (and hence $P\subset N_i'$ for some $i$), let $Q\in C_i$ be the point whose projection on $P_i$ is the center of $P$, and let $D_P$ be the spherical cap defined as the intersection of $S$ with the ball centered in $Q$ and of radius $\displaystyle \eps/2$. Similarly, define $D_{u(P)}$, corresponding to $u(P)$. Clearly, $D_P=\rho_P(D_{u(P)})$ for some $\rho_P\in SO(3)$. We remove from $M$ all the caps $D_{u(P)}$ and from $M'$ all the caps $D_P$, for $P\in \mathcal{P}'$.

Define now $s=\cala(M)$, $s'=\cala(M')$. Since $\sum {n_i'}\eps^2\rightarrow \sum \cala(N'_i)$, when $\eps\rightarrow 0$, we can choose $\eps$ and $1-m$ small enough such the above procedure removes from $M$ and $M'$ the sets $\mathcal{M}_1$ and $\mathcal{M}'_1$ of area greater than $\displaystyle \frac 12 s'$.

Inductively, define $S_i, S_i'$ as follows: $S_1=M\setminus \mathcal{M}_1$ and $S_1'=M'\setminus \mathcal{M}_1'$. By repeating the above process, obtain the sets $S_2, S_2'$ and so on.

Obviously, $\displaystyle \cala(S_t')<\left(\frac 12\right)^t\to0$ as $t$ grows to infinty. Since  $\cala(S_t)>s-s'>0$, there exists some $t$ such that $$\cala(S_t)>4\cdot\cala(S_t').$$ Once again, we go through the first step of the above construction applied to the sets $S_t$, $S_t'$ with the difference that $\mathcal{P}'$ will be the minimal set of all squares of lattices in $P_i$ which cover $\bigcup_{i=1}^kN_i'$, and $\mathcal{P}$ will contain all the squares of lattices with side length $2\eps$ that are included in $\bigcup_{i=1}^kN_i$. Also, $D_P$ will be the intersection of $S$ with the ball centered at $Q$ and of radius $\displaystyle\frac{\eps}{\sqrt{2}}$, and $D_{u(P)}$ is constructed analogously. The circle with the same center as $u(P)$, of radius $\displaystyle\frac{\eps}{\sqrt{2}}$, is included in $u(P)$.

Letting the set of $U_\al$ be the set of all $D_{u(P)}$, the conditions $(i)-(iii)$ in the Lemma are satisfied and the proof is complete.\qed


\begin{thebibliography}{99}
\bibitem{M} M. R. Murty, J. Esmonde, {\em Problems in algebraic number theory}, Springer-Verlag (2005)
\end{thebibliography}
\end{document}